\documentclass[11pt]{article}
\usepackage{amsmath}
\usepackage{mathrsfs}
\usepackage{amsfonts}

\usepackage{amsfonts, amsmath, amssymb}
\usepackage{amssymb,amsfonts,amsmath,color,
latexsym, epsfig,cite, psfrag,eepic,colordvi}
\usepackage{amscd,graphics}
\usepackage{lineno}
\usepackage{comment}

\newcommand{\qed}{\hfill $\Box $}

\newtheorem{theorem}{Theorem}[section]

\begin{document}

\title{A Degree Condition for  Graphs Having All $(a,b)$-Parity Factors}

\author{Haodong Liu and Hongliang Lu\footnote{luhongliang215@sina.com}\thanks{Supported by the National Natural
Science Foundation of China under grant No.11871391 and
Fundamental Research Funds for the Central Universities}\\ 
{School of Mathematics  and Statistics}
\\ {Xi'an Jiaotong University}\\
{Xi'an, Shaanxi 710049,  China}
}

\date{}

\maketitle

\begin{abstract}
Let $a$ and $b$ be positive integers such that $a\leq b$ and $a\equiv b\pmod 2$. 
We say that $G$ has all $(a, b)$-parity factors  if $G$ has an $h$-factor
 for every function $h: V(G) \rightarrow \{a,a+2,\ldots,b-2,b\}$ with $b|V(G)|$ even and $h(v)\equiv b\pmod 2$ for all $v\in V(G)$.
 In this paper, we prove that every graph $G$ with $n\geq 3(b+1)(a+b)$ vertices has all $(a,b)$-parity factors if $\delta(G)\geq (b^2-b)/a$, and for any two nonadjacent vertices $u,v \in V(G)$, $\max\{d_G(u),d_G(v)\}\geq \frac{bn}{a+b}$. Moreover, we show that this result is best possible in some sense.
\end{abstract}


\section{Introduction}

We consider simple graphs in this paper.
Let $G$ be a graph with vertex set $V(G)$ and edge set $E(G)$.
Given $x\in V (G)$, the set of vertices adjacent to $x$ is said to be the \emph{neighborhood} of $x$, denoted by $N_G(x)$, and $d_G(x)=|N_G(x)|$ is called the \emph{degree} of $x$. We write $N_G[x]$ for $N_G(x)\cup \{x\}$.
For a vertex set $D\subseteq V(G)$, let $N_D(v):=N_G(v)\cap D$ denote the set of vertices which are adjacent to $v$ in $D$.
For any vertex set $X\subseteq V(G)$, let $G[X]$ denote the vertex induced subgraph of $G$ induced by $X$ and the subgraph induced by vertex set $V(G)-X$ is also denoted by $G-X$.
A graph  $F$  is a \emph{spanning subgraph} of $G$ if $V(F)=V(G)$ and $E(F)\subseteq E(G)$. Let $h:V(G)\rightarrow \mathbb{N}$ and let  $S\subseteq V(G)$. We write $h(S):=\sum_{v\in S}h(v)$.

For a given graph $G$, let $g,f$ be two non-negative integer-valued functions such that $g(v)\leq f(v)$ and $g(v)\equiv f(v)\pmod 2$ for all $v\in V(G)$.
We say that a spanning subgraph $F$ of $G$ is a \emph{$(g,f)$-parity factor} if $d_F(v)\equiv f(v)\pmod 2$ and $g(v)\leq d_F(v)\leq f(v)$ for all $v\in V(G)$.
A $(g,f)$-parity factor is an \emph{$f$-factor} if $f(v)=g(v)$ for all $v\in V(G)$.
If $f(v)=k$ for all $v\in V(G)$, then an $f$-factor is a \emph{$k$-factor}.
Let $a,b$ be two integers such that $a\leq b$ and $a\equiv b\pmod 2$. If $f(v)=b$ and $g(v)=a$ for all $v\in V(G)$, then a $(g,f)$-parity factor is an \emph{$(a,b)$-parity factor}.
We call a graph $G$ \emph{having all $(a,b)$-parity factors} if $G$ has an $h$-factor for any functions $h$ with $a\leq h(v)\leq b$ and
$a\equiv h(v)\pmod 2$ for all $v\in V(G)$. If $G$ has an $h$-factor for any functions $h$ with $a\leq h(v)\leq b$ and $h(V(G))\equiv 0\pmod 2$, then we call $G$ having $(a,b)$-factors.

%
%
 Tutte \cite{Tutte} give a characterization for a graph to have an $f$-factor.
Lov\'asz \cite{Lov70,Lov72} gave a criterion for a graph to have an  $(g,f)$-parity factor. Kano and Tokushige \cite{KT92} gave a minimum degree condition for a graph to have all $(a,b)$-factors. In 1998, Niessen gave a  characterization for a graph to have \emph{all $(a,b)$-parity factors}.
\begin{theorem}[Tutte, \cite{Tutte}]\label{Tutte}
Let $G$ be a graph and let $f:V(G)\rightarrow N$.
$G$ contains an $f$-parity factor if and only if for all disjoint sets $S,T$ of $V(G)$,
\[
\eta(S,T)=f(S)-f(T)+d_{G-S}(T)-q_G(S,T;f)\geq 0,
\]
where $q_{G}(S,T;f)$ denotes the number of $f$-odd components $C$ of $G-S-T$ such that $f(V(C))+e_G(V(C),T)\equiv 1 \pmod 2$. Moreover, $\eta(S,T)\equiv f(V(G))\pmod 2$.
\end{theorem}
\begin{theorem}[Niessen, \cite{Niessen}]\label{Niessen}
Let $G$ be a graph and let $g,f:V(G)\rightarrow N$ such that $g(v)\leq f(v)$ and $g(v)\equiv f(v)\pmod 2$ for all $v\in V(G)$.
 $G$ has all $(g,f)$-parity factors if and only if
\[
\eta(S,T)=g(S)-f(T)+d_{G-S}(T)-q_G(S,T;f)\geq 0,
\]
for all disjoint sets $S,T$ of $V(G)$, where $q_{G}(S,T;f)$ denotes the number of $f$-odd components $C$ of $G-S-T$ such that $f(V(C))+e_G(V(C),T)\equiv 1 \pmod 2$.
\end{theorem}

Nishimura{\cite{Nishi}} gave a degree condition for a graph to have a $k$-factor.

\begin{theorem}[Nishimura, \cite{Nishi}]\label{Ni92}
Let $k$ be an integer such that $k\geq 3$, and let $G$ be a connected
graph of order $n$ with $n\geq 4k-3$, $kn$ even, and minimum degree at
least $k$. If $G$ satisfies
\begin{equation*}
\max\{d_G(u),d_G(v)\}\geq n/2
\end{equation*}
for each pair of nonadjacent vertices $u,v$ in $G$, then $G$ has a  $k$-factor.
\end{theorem}

%

In this paper, we give  a sufficient condition for a graph to have all $(a,b)$-parity factors.

\begin{theorem}\label{Main-Theorem}
Let $a,b,n$ be three integers such that $a\equiv b\pmod 2$, $a<b$, $na\equiv 0\pmod 2$  and $n\geq 3(b+1)(a+b)$.
Let $G$ be a connected graph of order $n$. If  $\delta(G)\geq \frac{b^2-b}{a}$ and
\begin{equation}\label{degree_condition}
max\{d_G(u),d_G(v)\}\geq \frac{bn}{a+b}
\end{equation}
for each pair of nonadjacent vertices $u$ and $v$ in $V(G)$, then $G$ has all $(a,b)$-parity factors.
\end{theorem}

\section{Proof of Theorem \ref{Main-Theorem}}

By contradiction, suppose that  the result does not hold.
By Theorem \ref{Niessen}, there exist two disjoint subsets $S$ and $T$ of $V(G)$ such that
\begin{align}\label{eq Nie1}
\eta(S,T)=a|S|-b|T|+\sum_{x\in T}d_{G-S}(x)-q(S,T)\leq -1,
\end{align}
where $q(S,T)=q(S,T;a)$ is the number of connected components $C$ of $G-S-T$ such that $a|V(C)|+e_G(V(C),T)\equiv 1 \pmod 2$ ($C$ is also called an $a$-odd components).
For simplicity, we write $s:=|S|$, $t:=|T|$ and $w:=q(S,T;a)$. By Theorem \ref{Tutte}, we have $\eta(S,T)\equiv a|V(G)|\equiv 0\pmod 2$. Hence (\ref{eq Nie1}) may be re-written as
\begin{align}\label{eq Nie2}
\eta(S,T)=as-bt+\sum_{x\in T}d_{G-S}(x)-w\leq -2.
\end{align}

We claim that  \begin{align}\label{ST-nonempty}
S\cup T\neq \emptyset.
\end{align} Otherwise, suppose that $s=t=0$. According to (\ref{eq Nie2}), we have $\eta(S,T)=-w\leq -2$, i.e.,  $w\geq 2$, a contradiction  since $G$ is connected.

When $w\geq 1$, let $C_1,C_2,\dots C_w$ denote the  $a$-odd components of $G-S-T$, and let   $m_i=|V(C_i)|$  for $1\leq i\leq w$.
Without loss of generality, suppose that $m_1\leq m_2\leq\ldots\leq m_w$.
Let $U=\bigcup_{1\leq i\leq w}V(C_i)$.

Without loss of generality, among all such subsets,  we choose $S$ and $T$ such that $U$ is minimal and $V(G)-S-T-U$ is maximal.

\medskip
\textbf{ Claim 1.~}   $d_{G-S}(u) \geq b+1$ and $ e_{G}(u,T) \leq a-1$ for each vertex $u\in U$.

By contradiction.
Firstly, suppose that there exists a vertex $u\in U$ such that $d_{G-S}(u) \leq b$. Let $T'=T\cup \{u\}$.
By inequality (\ref{eq Nie2}), we have
\begin{align*}
\eta(S,T')&=as-b|T'|+\sum_{x\in T'}d_{G-S}(x)-q(S,T')\\
&=(as-bt+\sum_{x\in T}d_{G-S}(x)-q(S,T'))+d_{G-S}(u)-b\\
&\leq as-bt+\sum_{x\in T}d_{G-S}(x)-q(S,T')\\
&\leq  as-bt+\sum_{x\in T}d_{G-S}(x)-(q(S,T)-1)\\
&\leq \eta(S,T)+1.
\end{align*}
Recall that  $\eta(S,T')\equiv \eta(S,T)\equiv 0\pmod 2$. So we have $\eta(S,T')\leq\eta(S,T)\leq -2$, which contradicts the choice of $U$.


Secondly, suppose that there exists a vertex $u\in U$ such that $e_{G}(u,T) \geq a$. Let $S'=S\cup \{u\}$. Then we have
\begin{align*}
\eta(S',T)&=a|S'|-bt+\sum_{x\in T}d_{G-S'}(x)-q(S',T)\\
          &=as+a-bt+\sum_{x\in T}d_{G-S}(x)-e_G(u,T)-q(S',T)\\
          &\leq as-bt+\sum_{x\in T}d_{G-S}(x)-q(S',T)\\
          &\leq as-bt+\sum_{x\in T}d_{G-S}(x)-(q(S,T)-1)\\
          &=\eta(S,T)+1.
\end{align*}
With similar discussion,  we have $\eta(S',T)\leq \eta(S,T)\leq -2$, which contradicts the choice of $U$. This completes the proof of Claim 1.

Note that when we choose appropriate $S$ and $T$, each components of $G-S-T$ contains at least $b-a+1$ vertices by Claim 1.

\medskip
\textbf{ Claim 2.~} Let $C_{i_1},\ldots,C_{i_{\tau}}$ be any $\tau$ components of $G[U]$ and let $U'=\bigcup_{j=1}^{\tau}V(C_{i_j})$.  $d_{G[T\cup U']}(u)\leq  b-1+\tau$ for every vertex  $u\in T$.

Suppose that there exists $u\in T$ such that $d_{G[T\cup U']}(u)\geq  b+\tau$. Let $T'=T-\{u\}$.
By (\ref{eq Nie2}), one can see that
\begin{align*}
\eta(S,T')&=as-b|T'|+\sum_{x\in T'}d_{G-S}(x)-q(S,T')\\
&= as-bt+b+\sum_{x\in T}d_{G-S}(x)-d_{G-S}(u)-q(S,T')\\
&\leq as-bt+b+\sum_{x\in T}d_{G-S}(x)-(b+\tau)-(q(S,T)-\tau)\\
&=as-bt+\sum_{x\in T}d_{G-S}(x)-q(S,T)\leq -2,
\end{align*}
contradicting to the maximality of $V(G)-S-T-U$. This completes the proof of Claim 2. \qed

By the definitions of $U$ and $C_i$, one can see that  $|U|\geq m_1+m_2(w-1)$ and so we have  $m_2\leq (|U|-m_1)/(w-1)$. Notice that $|U|+s+t\le n$. We have
\begin{equation}\label{m2}
m_2\leq \frac{n-s-t-m_1}{w-1}.
\end{equation}

By Claim 1, $e_{G}(u,T) \leq a-1$  for each $u\in U$. Thus for  $u\in U$, we have
\begin{equation}\label{eq-m_j-s-r}
  d_G(u)\leq (m_{j}-1)+s+r,
\end{equation}
where $r=\min\{a-1,t\}$.
Let $u_1\in V(C_1)$ and $u_2\in V(C_2)$. It follow from (\ref{eq-m_j-s-r}) that
\begin{equation}\label{eqm1m2}
    \max\{d_G(u_1),d_G(u_2)\}\leq (m_{2}-1)+s+r
\end{equation}
since $m_1\leq m_2$.


\medskip
\textbf{Claim 3.} $S\neq\emptyset$.

By contradiction. Suppose that $S=\emptyset$. Recall that $S\cup T\neq \emptyset$ by (\ref{ST-nonempty}). So we  have $t\geq 1$.
Since $b\geq a+2$ and $\delta(G)\geq (b^2-b)/a$, it follows by (\ref{eq Nie2}) that
\begin{align}\label{w>=(b+2)/a+2}
  w\geq \sum\limits_{v\in T}d_{G}(v)-bt+2\geq t(\frac{b^2-ab-b}{a})+2 \geq \frac{bt}{a}+2.
\end{align}

Combining (\ref{m2}) and (\ref{eqm1m2}), since $s=0$ and $r=\min\{a-1,t\}$, we have
\begin{align}\label{m2s-empty}
 \max\{d_G(u_1),d_G(u_2)\} \leq \frac{n-t-m_1}{w-1}+a-2.
\end{align}

By (\ref{w>=(b+2)/a+2}), we have
\[
\frac{n-t-m_1}{w-1}-1+a-1\leq \frac{a(n-t-m_1)}{bt+a}+a-2.
\]
  One can see that  $h(t)=\frac{n-t-m_1}{bt+a}$ is a monotone decreasing function on $t$. Since $t\ge 1$ and $n>\frac{(a-2)(b+a)}{b-a}$, we have \[\frac{n-t-m_1}{w-1}+a-2<\frac{an-a}{a+b}+a-2<\frac{bn}{a+b},\]
which implies that
\[
\max\{d_G(u_1),d_G(u_2)\}<\frac{bn}{a+b},
\]
contradicting  (\ref{degree_condition}).
This completes the proof of Claim 3.

\medskip
\textbf{Claim 4.} $T\neq\emptyset$.

Suppose that $T=\emptyset$, i.e., $t=0$.   
By Claim 1, $m_i\geq b+2$ for $1\leq i\leq w$.  By (\ref{eq Nie2}), we have $w\geq as+2.$ Thus
\[
n\geq w(b+2)+s\geq  (b+1)(as+2)+s,
\]
 i.e., $s \leq \frac{n-2(b+1)}{(b+1)a+1}$.
Recall that  $r=\min \{a-1,t\}=0$. We have
\begin{align*}
\max\{d_G(u_1),d_G(u_2)\}\leq m_2-1+s\leq \frac{n-s}{w-1}+s\leq \frac{n}{as+1}+s.
\end{align*}
 Let $h(s)=\frac{n}{as+1}+s$. Notice that $h(s)$ is a convex function on variable $s$. We know that the maximum value of the function $h(s)=\frac{n}{as+1}+s$ be obtained only at the boundary of $s$. Recall that $1\leq s< \frac{n}{a(b+1)+1}$. We have
\begin{align*}
h(s) &\leq \max\{h(1),h(\frac{n}{a(b+1)+1})\} \\
     &\leq \max\{\frac{n}{a+1}+1,\frac{n}{(b+1)a+1}+b+2\}.
\end{align*}
Notice that  $\frac{n}{a+1}+1<\frac{bn}{a+b}$ and $\frac{n}{(b+1)a}+b+2<\frac{bn}{a+b}$ since $n>3(a+b)$.
Hence we have
 $\max\{d_G(u_1),d_G(u_2)\}<\frac{bn}{a+b}$, which contradicts  (\ref{degree_condition}). \qed

Let $h_1:=\min\{d_{G-S}(v)\ |\ v\in T\}$, and let $x_1\in T$ be a vertex satisfying $d_{G- S}(x_1)=h_1$.
We write $p=|N_T[x_1]|$.
Furthermore, if $T- N_T[x_1]\neq\emptyset$, we put $h_2:=\min\{ d_{G- S}(v)\ |\ v\in T-N_T[x_1]\}$ and let $x_2\in T-N_T[x_1]$ such that $d_{G- S}(x_2)=h_2$.
Observe that $x_2$ and $x_1$ are not adjacent when $x_2$ exists and
\begin{align}\label{dx1x2}
 \max\{d_G(x_1),d_G(x_2)\}&\leq  \max\{h_1+s,h_2+s\}\leq h_2+s.
\end{align}


For completing the proof, we discuss four cases.

\medskip
\textbf{Case 1.} $h_1\geq b$.
\medskip

By (\ref{eq Nie2}), we have
\begin{align*}
  w &\geq as-bt+\sum_{v\in T}d_{G- S}(v)+2\\
 & \geq as+(h_1-b)t+2\nonumber\\
    & \geq as+2,\nonumber
\end{align*}
i.e.,
\begin{align}\label{case1-1}
  w &\geq  as+2.
\end{align}
Note that the number of vertices in graph $G$ satisfies $n\geq w+s+t$. Then from (\ref{case1-1}), we have
\begin{equation*}
  n \geq as+2+s+t
\end{equation*}
and it infers that
\begin{align*}
s<\frac{n-2}{a+1}.
\end{align*}
Recall that $h(s)=\frac{n}{as+1}+s$ is a convex function. Hence we have
\begin{align*}
\max\{d_G(x_1),d_G(x_2)\}&\leq m_2-1+s+r\\
&< \frac{n}{as+1}+s+a-2  \quad\mbox{(by Claim 1)}\\
&\leq \max\{h(1)+a-2,h(\frac{n-2}{a+1})+a-2\}\\ 
&\leq \max\{\frac{n}{a+1}+a-1,\frac{n-2}{a+1}+a\}\\
&<\frac{bn}{a+b}\quad \mbox{(since $n\geq a+b$)},
\end{align*}
which contradicts to (\ref{degree_condition}).

Now we may assume that $ h_1<b$ in the following discussion.

\medskip
\textbf{Case 2.} $T=N_T[x_1]$.
\medskip

Then we have $t=|N_T[x_1]|$ and $t\leq h_1+1 \leq b$. Note that  $\delta(G)\geq \frac{b^2-b}{a}$.  Since $s+h_1\geq \delta(G)$, we have $s\geq (\frac{b^2-b}{a}-h_1)$.
By (\ref{eq Nie2}), one can see that
\begin{align*}
  w&\geq as+(h_1-b)t+2\\
  &\geq as+(h_1-b)(h_1+1)+2\\
  &\geq a(\frac{b^2-b}{a}-h_1)+(h_1-b)(h_1+1)+2 \\
  &=(b-1)^2-(a+b-h_1-1)h_1+1\\
  &\geq (b-1)^2-(\frac{a+b}{2}-1)\frac{a+b}{2}+1 \quad \mbox{(since $a \equiv b \pmod 2$)}\\
 &\geq \frac{a+b}{2}+1 \quad \mbox{(since $a <b$)},
\end{align*}
i.e.,
\begin{align*}
  w \geq \frac{a+b}{2}+1.
\end{align*}
Now we have
\begin{align*}
\max\{d_G(u_1),d_G(u_2)\}&\leq m_2+s+a-1\\
&\leq \frac{n-s-t}{w-1}+s+a-2\\
&\leq \frac{2(n-s-t)}{a+b}+s+a-2.
\end{align*}
Let $h_1(s)= \frac{2(n-s)}{a+b}+s+a-2$. We discuss two subcases.

\medskip
\textbf{Case 2.1.~} $s\leq \frac{n}{3a}$.
\medskip

Notice that $h_1(s)$ is a linear function on $s$ and $1\leq s\leq n/3a$.  Then we have
\begin{align*}
h_1(s)&\leq h_1(\frac{n}{3a})\\
&=\frac{2(n-n/3a)}{a+b}+\frac{n}{3a}+a-2\\
&<\frac{bn}{a+b} \quad \mbox{(since $n\geq \frac{3(a+b)(a-2)}{2b-5}$)}.
\end{align*}
Thus we have
\begin{align*}
\max\{d_G(u_1),d_G(u_2)\} < \frac{bn}{a+b},
\end{align*}
a contradiction. 

\medskip
\textbf{Case 2.2.~} $s>\frac{n}{3a}$.
\medskip

By Claim 1, we have $m_1\geq b-a+1\geq 3$. Hence
\begin{align*}
n&\geq s+t+(b-a+1)w \\
&\geq s+t+3w\\
&\geq s+t+3as+3(h_1-b)t+6\quad \mbox{(by (\ref{eq Nie2}))}\\
&>n+\frac{n}{3a}-(3b-3h_1-1)t+6\\
&\geq n+\frac{n}{3a}-(3b-3h_1-1)(h_1+1)+6\quad \mbox{(since $t\leq h_1+1 \leq b$)}\\
   & \geq n+\frac{n}{3a}-3(\frac{b}{2}+\frac{1}{3})^2 \\
   & \geq n \quad\mbox{(since $n\geq a(\frac{3b}{2}+1)^2$)},
\end{align*}
 a contradiction.

Now  we may assume that $t>p$.

\medskip
\textbf{Case 3.~} $h_2\geq b$.
\medskip

Recall that $p=|N_T[x_1]|$. 
Notice that except $p$ vertices, each other vertex $v$ in $T$ satisfies $d_{G-S}(v)\ge h_2$. Thus by (\ref{eq Nie2}), we have
\begin{align}
w &\geq   as+\sum\limits_{v\in T}d_{G-S}(v)-bt+2\nonumber\\
&\geq as+(h_1-b)p+(h_2-b)(t-p)+2 \label{Ca3-w-b0}\\
  &\geq  as+(h_1-b)p+2\nonumber,
\end{align}
i.e.,
\begin{align}\label{Ca3-w-bound}
  w\geq  as+(h_1-b)p+2.
\end{align}

Now we discuss two subcases. 

\medskip
\textbf{Subcase 3.1.~} $h_2\leq \frac{(b+1)^2}{2}$.
\medskip

 From the assumption of theorem, we know that $\max\{d_G(x_1),d_G(x_2)\}\geq \frac{bn}{a+b}$, which implies that $h_2+s\geq  \frac{bn}{a+b}$. Then
\begin{equation}\label{subc31-eq1}
s\ge \frac{bn}{a+b}-h_2\geq \frac{bn}{a+b}-\frac{(b+1)^2}{2}.
\end{equation}
Now we have
\begin{align*}
n&\geq w+s+t\\
&\geq (a+1)s+(h_1-b)p+2+t\quad \mbox{(by (\ref{Ca3-w-bound}))}\\
&> (a+1)s+(h_1-b+1)p+2\quad \mbox{(since $p<t$)}\\
&\geq (a+1)s+(h_1-b+1)(h_1+1)+2\quad \mbox{(since $p\leq h_1+1\leq b$)}\\
&\geq (a+1)(\frac{bn}{a+b}-\frac{(b+1)^2}{2})+(h_1-b+1)(h_1+1)+2\quad \mbox{(by (\ref{subc31-eq1}))}\\
&\geq (a+1)(\frac{bn}{a+b}-\frac{(b+1)^2}{2})-\frac{b^2}{2}+2\\
&\geq n+2 \quad\mbox{(since $n\geq 2(b+1)(a+b)$)},
\end{align*}
a contradiction.

\medskip
\textbf{Subcase 3.2.~} $h_2>\frac{(b+1)^2}{2}$.
\medskip

Recall that $p\leq h_1+1\leq b$, $p<t$. Then we have
%
%
\begin{align*}
  w&\geq as+(h_1-b)p+(h_2-b)(t-p)+2 \quad \mbox{(by (\ref{Ca3-w-b0}))}\\
  &\geq as+(h_1-b)(h_1+1)+h_2-b+2\\
  &\geq as-\frac{(b+1)^2}{4}+\frac{(b+1)^2}{2}-b+2\\
  &\geq as+2,
\end{align*}
i.e.,
\begin{align}\label{subc32-eq1}
w\geq as+2.
\end{align}
By Claim 1 and (\ref{subc32-eq1}), we know that
\begin{align*}
n\geq s+t+(b-a+1)w\geq (3a+1)s+6,
\end{align*}
which implies that
\begin{align}\label{subc32-eq2}
s\leq \frac{n-6}{3a+1}.
\end{align}
Consider  $u_1\in V(C_1)$ and $u_2\in V(C_2)$. By (\ref{subc32-eq1}) and (\ref{subc32-eq2}), we have
\begin{align*}
\max\{d_G(u_1),d_G(u_2)\}&< m_2-1+r+s\\
&\leq \frac{n-s-t}{as+1}+s+a-2\\
&<\frac{n}{as+1}+s+a-2\\
&\leq \max\{\frac{n}{a+1}+a-1,\frac{n-6}{3a+1}+a+2\},
\end{align*}
where the last inequality holds since $1\leq s \leq \frac{n-6}{3a+1}$ and $h(s)=\frac{n}{as+1}+s$ is a convex function.
Observe that when $n\ge 2(a+b)$,  $\max\{\frac{n}{a+1}+a-1,\frac{n-6}{3a+1}+a+2\}<\frac{bn}{a+b}$. Thus we have
\begin{equation*}
\max\{d_G(u_1),d_G(u_2)\}< \frac{bn}{a+b}
\end{equation*}
when $n\ge 2(a+b)$, which contradicts  (\ref{degree_condition}).

\medskip
\textbf{Case 4.~} $0\leq  h_1\leq h_2\leq b-1$.
\medskip

By (\ref{degree_condition}) and (\ref{dx1x2}), since $x_1$ and $x_2$ are not adjacent, we have
$s\geq \frac{bn}{a+b}-h_2$.
By (\ref{eq Nie2}), we have
\begin{align}
  w &\geq  as+\sum\limits_{v\in T}d_{G-S}(v)-bt+2 \nonumber\\
  &\geq  as+(h_1-b)p+(h_2-b)(t-p)+2  \nonumber \\
  &= as+(h_2-b)t+(h_1-h_2)p+2.\label{eq:17}
\end{align}
One can see that
\begin{align}\label{case4}
  n&\geq s+t+w \nonumber \\
  &\geq (a+1)s+(h_2-b)t+(h_1-h_2)p+2+t \quad \mbox{(by (\ref{eq:17}))}\nonumber\\
  &= (a+1)s+(h_1-h_2)p+(h_2+1-b)t+2.
\end{align}
Note that $p\leq h_1+1$, $h_2+1-b\leq 0$, $s\geq \frac{bn}{a+b}-h_2$, $t\leq n-s\leq \frac{an}{a+b}+h_2$ and  $(h_2-h_1)(h_1+1)\leq h_2^2$. By taking these inequalities into  (\ref{case4}), we have
\begin{align*}
  n&\geq (a+1)(\frac{bn}{a+b}-h_2)+(h_1-h_2)(h_1+1)-(b-h_2-1)(\frac{an}{a+b}+h_2)+2\\
  &= (a+1)\frac{bn}{a+b}-(b-1)\frac{an}{a+b}+2+h_2(\frac{an}{a+b}+h_2-a-b)+(h_1-h_2)(h_1+1)\\
  &= n+2+h_2(\frac{an}{a+b}-a-b)+h_2^2-(h_2-h_1)(h_1+1)\quad \mbox{(since $n\geq \frac{(a+b)^2}{a}$)}\\
  &\geq n+2,
\end{align*}
a contradiction.


This completes the proof. \qed

\noindent\textbf{Remark.}
In Theorem \ref{Main-Theorem},  the bound in the assumption $\max\{d_{G}(u),d_{G}(v)\}\ge \frac{bn}{a+b}$ is tight.
For showing this, we construct following  graph:
Let $s,t,n$ be integers where $s=\lfloor \frac{bn-1}{a+b}\rfloor$ and $t=n-s$.
Let $K_n$ denote a complete graph of order $s$ and let $G=K_n-E(K_t)$.
Obviously, every pair of nonadjacent vertices $u,v$ in $G$ satisfies $max\{d_G(u),d_G(v)\} \geq \frac{bn}{a+b}-1$. Let $S=V(K_s)$ and $T=V(K_n)-S$.
Then we have
\[
\eta(S,T)=as-bt+\sum_{v\in T}d_{G-S}(v)-q(S,T)=as-bt<0,
\]
 i.e., $G$ doesn't satisfy the condition of Theorem \ref{Niessen}. Thus it is not true for graph $G$ to have all $(a,b)$-parity factors.


\end{document}